\documentclass{amsart}
\usepackage{lscape}
\newtheorem{theorem}{Theorem}[section]

\theoremstyle{definition}

\numberwithin{equation}{section}

\begin{document}
%======================================================================================================================================================
\title{Rank of elliptic curves associated to the Brahmagupta quadrilaterals}
%======================================================================================================================================================
\author{F. Izadi, F. Khoshnam, \and A. S. Zargar}
\address{Farzali Izadi: Department of Pure Mathematics, Faculty of Science, Urmia University, Urmia 165-57153, Iran.}
\email{f.izadi@urmia.ac.ir}
\address{Farzali Izadi: Department of Pure Mathematics, Azarbaijan Shahid Madani University, Tabriz 53751-71379, Iran.}
\email{farzali.izadi@azaruniv.edu}
\address{Foad Khoshnam: Department of Pure Mathematics, Azarbaijan Shahid Madani University, Tabriz 53751-71379, Iran.}
\email{khoshnam@azaruniv.edu}
\address{Arman Shamsi Zargar: Department of Pure Mathematics, Azarbaijan Shahid Madani University, Tabriz 53751-71379, Iran.}
\email{shzargar.arman@azaruniv.edu}
%======================================================================================================================================================
\subjclass[2000]{Primary 11G05, 14H52}
\keywords{Brahmagupta formula, Heron formula, quadrilaterals, Diophantine equation, elliptic curves, rank of elliptic curves}
%======================================================================================================================================================
\begin{abstract}
In this paper, we construct a family of elliptic curves of rank at least five. To do so, we use the  Brahmagupta formula for the area of cyclic quadrilaterals $(p^3,q^3,r^3,q^3)$ not necessarily standing for the genuine sides of quadrilaterals. It turns out that, as parameters of the curves, the integers $p,q,r,s$, along with the extra integers $u$, $v$ satisfy $u^6+v^6+p^6+q^6=2(r^6+s^6)$,  $uv=pq$. However, we utilize a subset of the solutions of the above system via the rational points of a specific elliptic curve of positive rank lying on the system.
\end{abstract}
%======================================================================================================================================================
\maketitle
%======================================================================================================================================================
\section{Introduction}
Let $E$ be an elliptic curve over $\Bbb Q$. The well-known theorem of Mordell and Weil states that $E(\Bbb Q)\simeq E(\Bbb Q)_{\rm tors}\times\Bbb Z^r$, where $r$ is a nonnegative integer called the Mordell-Weil rank of $E(\Bbb Q)$ and $E(\Bbb Q)_{\rm tors}$ is the subgroup of elements of infinite order called the torsion subgroup of $E(\Bbb Q)$. By a celebrated theorem of Mazur \cite{Sil1}, the only possible torsion subgroups over $\Bbb Q$ are $\Bbb Z/n\Bbb Z$ for $n=1,2,\ldots,10,12$ or $\Bbb Z/2\Bbb Z\times\Bbb Z/2n\Bbb Z$ for $1\leq n\leq4$. Even though there exists a folklore conjecture which says a rank can be arbitrarily high, it appears difficult to find examples of curves with high rank. The current record is an example of an elliptic curve over $\Bbb Q$ with rank $\geq28$ found by N. Elkies in May 2006, see \cite{D}. Besides, there is no known algorithm in order for determining the rank and also it is not known which integers can occur as ranks.

The notion of rank has attracted the interest of several authors in the last decades, and has led to some information on the structure of elliptic curves as $\Bbb Z$-modules. Specifically, there have been some investigations regarding ranks of certain elliptic curves via the Heron formula, see \cite{I1}.

In \cite{I1}, Izadi and Nabardi make a connection between the Heron formula for the area of triangles and elliptic curves. Recall that  the Heron formula for the area of a triangle $(A,B,C)$, say, is $$S=\sqrt{P(P-A)(P-B)(P-C)},$$
where $P$ is the semi-perimeter. Specifically, they use the Heron formula for the area of the formal triangle $(A^2,B^2,C^2)$ to generate infinitely many elliptic curves with high rank. More precisely, they show the elliptic curve
$$y^2=x^3+\frac{1}{4}\,(A^8+B^8+C^8-2A^4B^4-2A^4C^4-2B^4C^4)x$$
over the surface $A^4+D^4=2(B^4+C^4)$ has rank at least five and explicitly give the five independent points.

Throughout the paper, the curves we generate all have the trivial torsion subgroup $\mathcal T$. In this work, we use the Brahmagupta formula for the area of cyclic quadrilaterals to similarly find infinitely many elliptic curves with high rank. In effect, we consider the elliptic curve
$$
y^2=x^3-3(pqrs)^2x+2(pqrs)^3+\frac{1}{4}(p^6+s^6-q^6-r^6)^2-(p^3s^3+q^3r^3)^2,
$$
denoted by $E_{u,v,p,q,r,s}$, over 
\begin{align}
       &\quad u^6+v^6+p^6+q^6=2(r^6+s^6)\label{equation}\\
C:   & \nonumber\\
       & \quad uv=pq \label{eq}
\end{align}
and prove that the group of the rational map $C\rightarrow E_{u,v,p,q,r,s}$, that commutes with the projection $E_{u,v,p,q,r,s}\rightarrow C$, has rank at least five. We do this by exhibiting five explicit sections $P_1$, $P_2$, $P_3$, $P_4$, $P_5$, and showing that these are linearly independent. We use the fact that $C(\Bbb Q)$, the set of rational solutions satisfying on $C$, is infinite in order to deduce that infinitely many specializations of $E_{u,v,p,q,r,s}$ have ranks at least five over the rationals. This is done by using $$y^2-28xy-560y=x^3-20x^2-400x+8000,$$
an elliptic curve of positive rank lying on $C$ found in \cite{I2}. In fact, we show the following theorem.
\begin{theorem}\label{main}
There are infinitely many elliptic curves over $C$ of rank at least five, parameterized by an elliptic curve of rank at least three over $\Bbb Q(p,q,r,s)$.
\end{theorem}

\section{The construction of $E_{u,v,p,q,r,s}$}
As is mentioned above, in this work we deal with elliptic curves related to (positive) integer solutions on $C$. The Brahmagupta formula for the area of a cyclic quadrilateral in terms of its side lengths states that for a quadrilateral with sides $a,b,c,d$, the area of the quadrilateral is given by
\begin{equation}\label{hq1}
S=\sqrt{(P-a)(P-b)(P-c)(P-d)},
\end{equation}
in which $P=(a+b+c+d)/2$ is the semi-perimeter. Equivalently,
\begin{equation}\label{hq2}
16S^2=(b+c+d-a)(a+c+d-b)(a+b+d-c)(a+b+c-d).
\end{equation}
By an appropriate arrangement, the formula \eqref{hq2} becomes
\begin{equation}\label{hq3}
16S^2=\left[(a+d)^2-(c-b)^2\right]\left[(b+c)^2-(d-a)^2\right].
\end{equation}
By some algebra, \eqref{hq3} turns out to be
$$
16S^2=(2ad+2bc)^2-(a^2+d^2-b^2-c^2)^2,
$$
or
\begin{equation}\label{hq4}
\frac{1}{4}\left(b^2+c^2-a^2-d^2\right)^2=(ad+bc)^2-4S^2.
\end{equation}
Now, take $(a,b,c,d)=(p^3,q^3,r^3,s^3)$. By substituting, \eqref{hq4} turns into
\begin{equation}\label{hq5}
\frac{1}{4}\left(q^6+r^6-p^6-s^6\right)^2=(p^3s^3+q^3r^3)^2-4S^2.
\end{equation}
Expanding and sorting the right hand side of \eqref{hq5}, we get
$$\frac{1}{4}\left(q^6+r^6-p^6-s^6\right)^2=(p^2s^2)^3+(q^2r^2)^3+2(pqrs)^3-4S^2,$$
or
\begin{equation}\label{hq6}
\left(\frac{q^6+r^6-p^6-s^6}{2}\right)^2=(p^2s^2+q^2r^2)^3-3p^2q^2r^2s^2(p^2s^2+q^2r^2)+2(pqrs)^3-4S^2.
\end{equation}
Setting $x=p^2s^2+q^2r^2$, $y=(q^6+r^6-p^6-s^6)/2$, and using \eqref{hq5}, we can define the following elliptic curve
\begin{equation}\label{curve0}
y^2=x^3-3(pqrs)^2x+2(pqrs)^3+\frac{1}{4}(p^6+s^6-q^6-r^6)^2-(p^3s^3+q^3r^3)^2,
\end{equation}
denoted by $E_{p,q,r,s}$, over $\Bbb Q(p,q,r,s)$.

By symmetric roles of $p, q, r, s$ in our formulas, the elliptic curve has the following non-obvious points:
\begin{align*}
P_1(p,q,r,s)&=\left(p^2q^2+r^2s^2,\frac{p^6+q^6-r^6-s^6}{2}\right), \\
P_2(p,q,r,s)&=\left(p^2r^2+q^2s^2,\frac{p^6-q^6+r^6-s^6}{2}\right), \\
P_3(p,q,r,s)&=\left(p^2s^2+q^2r^2,\frac{p^6-q^6-r^6+s^6}{2}\right).
\end{align*}

By the specialization theorem (\cite{Sil2} or \cite[Theorem 11.4, p. 271]{Sil1}), in order to prove that the family of elliptic curves defined in \eqref{curve0} has rank at least three over $\Bbb Q(p,q,r,s)$, it suffices to find a specialization $p=p_0$, $q=q_0$, $r=r_0$, $s=s_0$ such that the points $P_i(p,q,r,s)$, $i=1,2,3$, are linearly independent on the specialized curve over $\Bbb Q$. If we take $p=1$, $q=1/2$, $r=2$, $s=0$, then the points $P_1(1,1/2,2,0)=\left(1/4,4031/128\right)$, $P_2(1,1/2,2,0)=\left(4,4159/128\right)$, $P_3(1,1/2,2,0)=\left(1,4033/128\right)$ are linearly independent points of infinite order on the elliptic curve
$$E_{1,1/2,2,0}:y^2=x^3+16248705/16384.$$
Indeed, the determinant of the N\'{e}ron-Tate height pairing matrix of these points is the nonzero value $170.021501512688$, according to SAGE \cite{Sag}. We note that rank $E_{1,1/2,2,0}(\Bbb Q)=4$ carried out by Cremona's  \verb+mwrank+ program \cite{C}.

\section{Infinitude of curves $E_{u,v,p,q,r,s}$ with rank at least five}
In order to increase the rank of \eqref{curve0}, we set $u$ and $v$ satisfy \eqref{equation}--\eqref{eq}, i.e.,
$$u^6+v^6+p^6+q^6=2(r^6+s^6), \ uv=pq.$$
Note that since the quadrilateral $(p,q,r,s)$ comes from \eqref{equation}--\eqref{eq}, there is no guarantee that $(p,q,r,s)$ and $(p^3,q^3,r^3,s^3)$ are genuine ones.

In \cite{I2}, the authors show the existence of infinitely many integer solutions to \eqref{equation}. The method is based on the points of the elliptic curve
$$E:y^2-28xy-560y=x^3-20x^2-400x+8000,$$
being generated by $(-50,400)$, explained as follows: Suppose $G=(A/B^2,C/B^3)$ is a rational point on the elliptic curve $E$ with $A,B,C \in \Bbb Z$ and $\gcd(B,AC)=1$. Then,
$$\left\{\begin{array}{l}
u=10B(A-20B^2)(-C+4BA-80B^3) \\
\qquad\times(-C^2+8ACB-160CB^3-20A^2B^2-400AB^4+8000B^6+A^3), \vspace{.2cm} \\
v=C(10BA-200B^3-C) \\
\qquad\times(-C^2+8ACB-160CB^3-20A^2B^2-400AB^4+8000B^6+A^3), \vspace{.2cm} \\
p=C(-2880000B^9+5360B^4CA-C^3+26AC^2B-164B^2CA^2+2BA^4 \\
\qquad-49600B^6C+CA^3-520C^2B^3+240B^3A^3-19200B^5A^2+416000B^7A), \vspace{.2cm} \\
q=10B(A-20B^2)(-C+4BA-80B^3) \\
\qquad\times (C^2-20A^2B^2-400AB^4+8000B^6+A^3), \vspace{.2cm} \\
r=(10BA-200B^3-C)(-C+4BA-80B^3) \\
\qquad\times(-C^2+8ACB-160CB^3-20A^2B^2-400AB^4+8000B^6+A^3), \vspace{.2cm} \\
s=6BC(A-20B^2)\\
\qquad\times(-C^2+8ACB-160CB^3-20A^2B^2-400AB^4+8000B^6+A^3).
\end{array}\right.$$
is an integral solution to \eqref{equation}.

One can readily observe that
\begin{align*}
uv-pq&=20BC(A-20B^2)^4(-C+4BA-80B^3)^2\\
&\qquad \times(A^3-20B^2A^2-400B^4A+28CBA+8000B^6+560CB^3-C^2).
\end{align*}
But the point $(A/B^2,C/B^3)$ lies on $E$, i.e.,
$$A^3-20B^2A^2-400B^4A+28CBA+8000B^6+560CB^3-C^2=0.$$
Hence, $uv=pq$.

From now on, the curve \eqref{curve0} is denoted by $E_{u,v,p,q,r,s}$. We first show that the following points
\begin{align*}
P_4(u,v,p,q,r,s)&=\left(u^2r^2+v^2s^2,\frac{u^6-v^6+r^6-s^6}{2}\right), \\
P_5(u,v,p,q,r,s)&=\left(u^2s^2+v^2r^2,\frac{u^6-v^6-r^6+s^6}{2}\right),
\end{align*}
lie on $E_{u,v,p,q,r,s}$: The point $P_4(u,v,p,q,r,s)$ satisfies $E_{u,v,p,q,r,s}$ if and only if we have
\begin{align*}
&(u^{12}+v^{12}-p^{12}-q^{12})\\
&\quad -2(v^6s^6-p^6s^6-p^6q^6-p^6r^6+u^6r^6-s^6q^6+u^6v^6+u^6s^6-q^6r^6+r^6v^6)\\
&\quad -12(u^4r^4v^2s^2+u^2r^2v^4s^4-p^2q^2r^4s^2u^2-p^2q^2r^2s^4v^2)=0,
\end{align*}
or, equivalently, by \eqref{eq},
\begin{align*}
&(q^2-u^2)(q^2+qu+u^2)(q^2-qu+u^2)(p^2-u^2)(p^2+pu+u^2)(p^2-pu+u^2)\\
&\quad \times (p^6q^6+p^6u^6-2u^6r^6-2s^6u^6+q^6u^6+u^{12})=0.
\end{align*}
Brushing aside the impossible cases, it is equivalent to having
$$p^6q^6+p^6u^6-2u^6r^6-2s^6u^6+q^6u^6+u^{12}=0,$$
or, equivalently
$$u^{12}+(p^6-2r^6-2s^6+q^6)u^6+p^6q^6=0,$$
which by using \eqref{equation}, this is the same
$$u^{12}-(u^6+v^6)u^6+p^6q^6=0,$$
or, equivalently, $u^6v^6-p^6q^6=0$ which holds by \eqref{eq}. Similarly, one can see $P_5(u,v,p,q,r,s)$ lies on $E_{u,v,p,q,r,s}$.

Now, specialization at $(u,v,p,q,r,s)=(748,6825,6188,825,6545,468)$ shows these five points are independent.
In fact for this value of $(u,v,p,q,r,s)$, the specialized curve $E_{748,6825,6188,825,6545,468}$, i.e.,
\begin{align*}
y^2& = x^3 - 733568661605545473708000000x \\
&\qquad+
126120106752478587342832201925322189643214400
\end{align*}
with the points
\begin{align*}
P_1(748,6825,6188,825,6545,468)&=(35444382573600, 11231145747112000830120), \\
P_2(748,6825,6188,825,6545,468)&=(1640436333421600, 67374938649066874419880), \\
P_3(748,6825,6188,825,6545,468)&=(37542673468881, 11231450539996549524921), \\
P_4(748,6825,6188,825,6545,468)&=(34169761645600, 11230981104248186419880), \\
P_5(748,6825,6188,825,6545,468)&=(1995497942444721, 89837370293311610364681)
\end{align*}
has nonzero regulator $1.29010146713893\times 10^7$, and we have
$$E_{748,6825,6188,825,6545,468}\simeq\Bbb Z^5.$$
The concept of the specialization theorem, thus, shows the existence of infinitely many elliptic curves of the form $E_{u,v,p,q,r,s}$ of rank at least five over $u^6+v^6+p^6+q^6=2(r^6+s^6)$ and trivial torsion subgroup $\mathcal T$.

\bibliographystyle{amsplain}

\begin{thebibliography}{10}

\bibitem{C} J. Cremona, http://maths.nottingham.ac.uk/personal/jec/ftp/progs.

\bibitem{D} A. Dujella, http://web.math.pmf.unizg.hr/\verb+~+duje/tors.html

\bibitem{I1} F. Izadi, K. Nabardi, \textit{A family of elliptic curves with rank $\geq5$}, to appear in Period. Math. Hungar.

\bibitem{I2} F. Izadi, F. Khoshnam, A. S. Zargar, \textit{A note on the Diophantine equations $x_1^k+x_2^k+x_3^k+x_4^k=2(y_1^k+y_2^k), k=3, 6$}, Notes Number Theory Discrete Math., {\bf 20} (2014), 1--10.

\bibitem{Sag} Sage software, version 4.5.3, available at http://www.sagemath.org.

\bibitem{Sil1} J. H. Silverman, \textit{Advanced Topics in the Arithmetic of Elliptic Curves}, Springer, New York, 1994.

\bibitem{Sil2} J. H. Silverman, \textit{Heights and the specialization map for families of abelian varieties}, J. Reine Angew. Math., {\bf 342} (1983), 197--211.
\end{thebibliography}

\end{document}